\documentclass[12pt,leqno]{article}
\usepackage{amsmath,amsthm}

\def\init{\setcounter{equation}{0}}
\newtheorem{theorem}{Theorem}[section]

\newcommand{\R}{{\bf R}}

\newtheorem{lemma}{Lemma}[section]

\newcommand{\e}{{\varepsilon}}

\title{Uniqueness and nonuniqueness in inverse hyperbolic problems and the black holes
phenomenon.
\author{G.Eskin, \ \ \  Department of Mathematics, UCLA,\\ Los Angeles,
CA 90095-1555, USA. \ E-mail: eskin@math.ucla.edu}
}

\begin{document}

\maketitle

\begin{abstract}
This paper consists of two parts.  In the first part
  we describe
the recent works   
 on the inverse problems 
for the wave equation in $(n+1)$-dimensional 
space equipped with pseudo-Riemannian metric with Lorentz signature.
We study the conditions of the existence of black (or white) holes for these 
wave equations.  
In the second part we prove energy type estimates on a finite time interval in the presence
of black or white holes.  We use these estimates to prove the nonuniqueness of the inverse problems.
\end{abstract}

\section{Introduction.}
\label{section 1}
\init
 
A powerful method for solving the inverse hyperbolic problem for equations 
of the form $\frac{\partial^2u}{\partial t^2}+Au=0$  where $A$ is a Laplace-Beltrani operator 
with time-independent coefficients, was discovered by M.Belishev more then twenty years ago.  
It is called the Boundary Control (BC) method.
It was further developed by M.Belishev,  M.Belishev and Y.Kurylev, Y.Kurylev and M.Lassas 
and others  (see [B1], [B2], [KKL] and further references there).
An impotant part of the solution of the hyperbolic inverse problem is played by the unique
continuation theorem due to D.Tataru [T].
In [E1], [E2] the author proposed a new approach to the inverse hyperbolic problem that
includes  ideas from the BC-method.  This new
method allowed to solve some inverse hyperbolic problems that were not accessible by
the BC-method:  In [E3] the case of hyperbolic equations with time-dependent coefficients was 
considered and in [E4] the case of the hyperbolic equation with general pseudo-Riemannian 
time-independent metric was treated.  In the following sections we describe the main results of
[E4]  and [E5].

An interesting phenomenon discussed in [E5]  is the appearance of black holes.  
These black holes are  called artificial black holes (they are also called acoustic black holes,
or optical black holes)  to distinguish
from the black holes in the general relativity.
Artificial black holes attracted a great interest of physicists (see [NVV], [V] and
additional references there)  because the physisists hope to create and study the black hole in
the laboratory and expect that this will help in the understanding of the black holes in
the universe.

In the last two sections we prove the energy type estimates on a finite time interval in the presence
of black or white holes.  We use these estimates to  prove the nonuniqueness 
in the inverse  problems.

\section{The inverse hyperbolic problems.}
\label{section 2}
\init

Let $\Omega$ be a bounded domain in $\R^n$  with smooth boundary $\partial\Omega$.  Let 
$\Gamma\subset\partial\Omega$  be  an open subset of $\partial\Omega$.

Consider a hyperbolic equation in the cylinder $\Omega\times\R$:
\begin{equation}                                               \label{eq:2.1}
\sum_{j,k=0}^n\frac{1}{\sqrt{(-1)^ng(x)}}\frac{\partial}{\partial x_j}\left(\sqrt{(-1)^ng(x)}g^{jk}(x)
\frac{\partial u(x_0,x)}{\partial x_k}\right)=0,
\end{equation}
where
$x=(x_1,...,x_n)\in \overline{\Omega},\ x_0\in \R$ is the time variable, 
the coefficients in (\ref{eq:2.1}) are smooth and independent of $x_0$,
$[g_{jk}(x)]_{j,k=0}^n=([g^{jk}(x)]_{j,k=0}^n)^{-1}$
is a pseudo-Riemannian metric with the Lorenz signature,  i.e. the quadratic form
$\sum_{j,k=0}^ng^{jk}(x)\xi_j\xi_k$
has the signature $(1,-1,-1,...,-1)$  for all $x\in \overline{\Omega}$.
Here $g(x)=(\det[g^{jk}(x)]_{j,k=0}^n)^{-1}$.
Note that $(-1)^ng(x)>0,\ \forall x\in \overline{\Omega}$.

  We assume  that
\begin{equation}                                               \label{eq:2.2}
g^{00}(x)>0,\ \ x\in \overline{\Omega},
\end{equation}
i.e. $(1,0,...,0)$ is not a characteristic direction,
and that 
\begin{equation}                                              \label{eq:2.3}
\sum_{j,k=1}^ng^{jk}(x)\xi_j\xi_k<0\ \ \ \ \mbox{for\ \ }\forall (\xi_1,...,\xi_n) \neq
(0,...,0),\ \forall x\in \overline{\Omega},
\end{equation}
i.e.  the quadratic form  (\ref{eq:2.3})   is negative definite.
Note that (\ref{eq:2.3}) equivalent to the condition that
\begin{equation}                                               \label{eq:2.4}
g_{00}(x)>0,\ \ x\in\overline{\Omega},
\end{equation}
i.e. that $(1,0,...,0)$  is a time-like direction.

We consider the initial-boundary value problem for the equation (\ref{eq:2.1}) in the
cylinder $\Omega\times \R$:
\begin{equation}                                            \label{eq:2.5}
u(x_0,x)=0\ \ \ \mbox{for}\ \ \ x\in\Omega,\ x_0<< 0,
\end{equation}
\begin{equation}                                         \label{eq:2.6}
u(x_0,x)\left|_{\partial\Omega\times\R}\right. = f(x_0,x'),\ x'\in\partial\Omega,
\end{equation}
where
$f(x_0,x')$  has a compact support in $\partial\Omega\times\R$.

Let $\Lambda f$  be the Dirichlet-toNeumann (DN) operator, i.e.
\begin{equation}                                            \label{eq:2.7}
\Lambda f=\sum_{j,k=0}^ng^{jk}(x)\frac{\partial u}{\partial x_j}\nu_k(x)
\left.(-\sum_{p,r=1}^ng^{jk}(x)\nu_p\nu_r)^{-\frac{1}{2}}\right|_{\partial\Omega\times\R},
\end{equation}
where
$\nu_0=0,(\nu_1,...,\nu_n)$  is the unit outward normal vector to $\partial\Omega\subset \R^n$,  
$u(x_0,x)$  is the solution of 
(\ref{eq:2.1}), (\ref{eq:2.5}),  (\ref{eq:2.6}).

Consider a smooth change of variables of the form:
\begin{eqnarray}                              \label{eq:2.8}
\hat{x}_0=x_0+a(x),
\\
\nonumber
\hat{x}=\varphi(x),
\end{eqnarray}
where $\varphi(x)$  is a diffeomorphism of $\overline{\Omega}$
onto some domain $\overline{\hat{\Omega}}$ such that $\overline{\Gamma}\subset\partial\hat{\Omega}$,
$\varphi(x)=x$  on $\overline{\Gamma}$, $a(x)=0$  on $\overline{\Gamma}$.   Note that (\ref{eq:2.8}) is
an identity map on $\overline{\Gamma}\times \R$.
Note also that the map (\ref{eq:2.8}) transforms (\ref{eq:2.1}) into an equation of the same form in
$\hat{\Omega}\times\R$.

The following theorem holds:
\begin{theorem} (c.f.[E4]):                                  \label{theo:2.1}
Let $L$ and $\hat{L}$  be two operators of the form (\ref{eq:2.1})  in $\Omega\times \R$
and $\hat{\Omega}\times\R$  respectively.
Consider initial-boundary value problems of the form (\ref{eq:2.5}), (\ref{eq:2.6})  for
$L$ and $\hat{L}$.  Suppose $\Lambda f=\hat{\Lambda}f$  on $\Gamma \times \R$  for all 
$f\in C_0^\infty(\Gamma\times\R)$  where $\Lambda, \hat\Lambda$ are DN operators for $L,\hat L$,
respectively.  Suppose that conditions (\ref{eq:2.2})  and (\ref{eq:2.3}) hold for $L$ and $\hat L$.
Then there exists a map of the form (\ref{eq:2.8}) such that 
\begin{equation}                                       \label{eq:2.9}
[\hat g^{jk}(\hat x)]_{j,k=0}^n=J^T(x)[g^{jk}(x)]_{j,k=0}^nJ(x),
\end{equation}
where 
$([\hat g^{jk}(\hat x)]_{j,k=0}^n)^{-1}$
is the metric tensor for $\hat L$  and $J(x)$  is the Jacobi matrix of (\ref{eq:2.8}).
\end{theorem}

{\bf Remark 2.1.}
It is enough to know the DN operator on $\Gamma\times(0,T_0)$  for some $T_0>0$  instead of
$\Gamma\times\R$.
More precisely,  let $T_+$  be the smallest number such that
$D_+(\overline{\Gamma}\times\{x_0=0\})\supset\overline{\Omega}\times\{x_0=T_+\}$  where 
$D_+(\overline{\Gamma}\times\{x_0=0\})$  is the forward domain of influence of
$\overline{\Gamma}\times\{x_0=0\}$
corresponding to (\ref{eq:2.1}).
Analogously let $T_-$  be the smallest number such that
$D_-(\overline{\Gamma}\times\{x_0=T_-\})\supset\overline{\Omega}\times\{x_0=0\}$  where  
$D_-(\overline{\Gamma}\times\{x_0=T_-\})$  
is the backward domain of influence of $\overline{\Gamma}\times\{x_0=T_-\}$.
If $T_0>T_-+T_+$  then $\Lambda=\hat\Lambda$  on $\Gamma\times(0,T_0)$  implies (\ref{eq:2.9}),
i.e.  the isometry of metrics $[g_{jk}(x)]$  and $[\hat g_{jk}(\hat x)]$.

\section{The equation of the propagation of light in the moving dielectric medium.}
\label{section 3}
\init

In this section we apply Theorem \ref{theo:2.1} to the equation of the propagation of light 
in the moving medium.

It was discovered by Gordon [G]  that the equation of the propagation of light in
a moving medium is given by the hyperbolic equation  of the form (\ref{eq:2.1})  with the metric tensor
\begin{equation}                                      \label{eq:3.1}
g^{jk}(x)=\eta^{jk}+(n^2(x)-1)v^j(x)v^k(x),
\end{equation}
$0\leq j,k \leq n,\ n=3$,  where
$[\eta_{jk}]=[\eta^{jk}]^{-1}$  is the Lorentz metric tensor:
$\eta^{jk}=0$ when $j\neq k,\ \eta^{00}=1,\ \eta^{jj}=-1,\ 1\leq j\leq n,\ x_0=t$
is the time,  
$n(x)=\sqrt{\e(x)\mu(x)}$  is the refraction index,
$w(x)=(w_1(x),w_2(x),w_3(x))$   is the velocity of
flow,
$$
v^{(0)}=
\left(1-\frac{|w|^2}{c^2}\right)^{-\frac{1}{2}},\ \ \
v^{(j)}=
\left(1-\frac{|w|^2}{c^2}\right)^{-\frac{1}{2}}\frac{w_j(x)}{c},\ \ 1\leq j\leq 3,
$$
  is the  four-velocity field 
 of the flow,
 $c$ is the speed of light in the vacuum. 
We shall call the equation (\ref{eq:2.1}) with metric (\ref{eq:3.1})
the Gordon equation.

Let
$\Omega$ be a smooth domain in $\R^n$ of the form  
$\Omega=\Omega_0\setminus\cup_{j=1}^m\overline{\Omega}_j$  where $\Omega_0$  is 
simply-connected,  $\Omega_j,\ 1\leq j\leq m,$  are smooth domains called obstacles,
$\overline{\Omega}_j\subset\Omega_0,\ 1\leq j\leq m,$
$\overline{\Omega}_j\cap\overline{\Omega}_k=\emptyset$  when 
$j\neq k$.

We shall consider the following initial-boundary value problem  for the Gordon equation:
\begin{eqnarray}                                       \label{eq:3.2}
u(x_0,x)=0\ \ \ \mbox{for}\ \ x_0<<0,\ x\in\Omega,
\\
\nonumber
u(x_0,x)|_{\partial\Omega_j\times\R}=0,\ \ 1\leq j\leq\ m,
\\
\nonumber                                       
u(x_0,x)|_{\partial\Omega_0\times\R}=f(x_0,x),
\end{eqnarray}
i.e.
$\partial\Omega_0=\Gamma$ in the notations of Theorem \ref{theo:2.1}.

Note that the condition (\ref{eq:2.2}) of \S 2 is always satisfied  since
$$
g^{00}=1+(n^2-1)(v^0)^2>0.
$$

The condition that any direction $(0,\xi_1,...,\xi_n)$  is not characteristic (c.f. 
condition (\ref{eq:2.3})) holds when
\begin{equation}                                 \label{eq:3.3}
|w(x)|^2<\frac{c^2}{n^2(x)}.
\end{equation}
We shall impose some restrictions on the flow $w(x)$.  Let $x=x(s)$  be a trajectory of the
flow,  i.e.
$$
\frac{dx}{ds}=w(x(s)),\ \ \ 0\leq s\leq 1,
$$
where $w(x(s))\neq 0$  for $0\leq s\leq 1$.  
We assume the following condition:
\begin{eqnarray}               
\nonumber
\mbox{(A)\ \ \ 
 The trajectories that start and end
on}\   \partial\Omega_0,\ \  
\\
\nonumber
\ \ \ \ \ \ \ \ \   \mbox{or are closed curves in}\ \  \Omega,
\ 
\mbox{are dense in}\ \  \overline{\Omega}.
\end{eqnarray}

\begin{theorem}[(c.f. [E5])]                                   \label{theo:3.1}
Let $[g_{jk}(x)]_{j,k=0}^n$ and $[\hat{g}_{jk}(y)]_{j,k=0}^n$
be two Gordon metrics in domains $\Omega$  and $\hat{\Omega}$,  respectively.
Consider two initial-boundary value problems of the form (\ref{eq:3.2}) in
$\Omega\times\R$  and $\hat{\Omega}\times\R$,  respectively,  
where $\Omega=\Omega_0\setminus\cup_{j=1}^m\overline{\Omega}_j$,
$\hat\Omega=\Omega_0\setminus\cup_{j=1}^{\hat m}\overline{\hat\Omega}_j$.  
Assume that 
the refraction indexes  $n$  and $\hat{n}$ are constant and  that
the flow $w(x)$ satisfies the condition (A). 
Assume also that (\ref{eq:3.3}) holds for both metrics.  Then
$\Lambda=\hat\Lambda$  on $\partial\Omega_0\times\R$  implies that 
$\hat n=n,\ \hat\Omega=\Omega$  and  the flows $\hat w(x)$  and $w(x)$ are equal.
\end{theorem}

\section{The propagation of light in the slowly moving medium.}
\label{section 4}
\init

In this case one drops the terms of order  $\frac{|w|^2}{c^2}$.  Then the metric tensor has the form:
\begin{eqnarray}                                             \label{eq:4.1}
g^{jk}(x)=\eta^{jk}\ \ \mbox{for}\ \ 1\leq j,k\leq n,\ n=3,                
\\
\nonumber
g^{00}(x)=n^2(x),\ \ g^{0j}(x)=g^{j0}(x)=(n^2(x)-1)\frac{w_j(x)}{c},\ \ 1\leq j\leq n.
\end{eqnarray}
The wave equation with metric (\ref{eq:4.1}) describes the propagation of light in 
 a slowly moving medium.  
We shall see that the inverse problem for such equation exhibits some nonuniqueness.

Denote $v_j(x)=g^{0j}=g^{j0}$.  We say that the flow $v=(v_1,...,v_n)$  is a gradient flow
if $v(x)=\frac{\partial b(x)}{\partial x}$  where $b(x)\in C^\infty(\overline{\Omega}),\ b(x)=0$
on $\partial\Omega_0$.

\begin{theorem}(c.f. [E4])                                         \label{theo:4.1}
Consider two initial-boundary value problems in domains $\Omega\times\R$  and 
$\hat\Omega\times \R$  for operators of the form (\ref{eq:2.1})  with metrics $[g^{jk}(x)],
\ [\hat g^{jk}(\hat x)]$ of the form (\ref{eq:4.1}).  Assume that the DN operators $\Lambda$
and $\hat\Lambda$  are equal on $\partial\Omega_0\times \R$.  Assume that there exists an open
connected and dense $O\subset\Omega$  such that $v(x)$ does not vanish on $O$.  Then
$\hat\Omega=\Omega,\ \hat n(x)=n(x)$  and $\hat{v}(x)=v(x)$  if $v(x)$  is not a gradient flow.
In the case of the gradient flow there are two solutions of the inverse problem:
$$
\hat v(x)=v(x)\ \ \ \mbox{and}\ \ \hat v(x)=-v(x).
$$
\end{theorem}

{\bf Remark 4.1} (c.f. [E4])
Suppose that the open set $O$ where $v(x)\neq 0$ consists of several open components
$O_1,...,O_r$.   Suppose there exists $b_j(x)\in C^\infty(\overline{\Omega}),\ 
b_j(x)=0$  on $\partial\Omega_0,\ \frac{\partial b_j}{\partial x}=v(x)$ on $O_j,\ b_j=0$
in $\overline{\Omega}\setminus O_j,\ j=1,2,...,r$.

Then we have $2^r$  solutions of the inverse problem where each of these solutions is 
equal to either $\frac{\partial b_j}{\partial x}$  or to $-\frac{\partial b_j}{\partial x}$
on $O_j$.  

\section{Artificial black holes.}
\label{section 5}
\init

Let $S(x)=0$  be a smooth closed surface in $\R^n$ such that the surface $S\times\R\subset\R^{n+1}$  is
 a characteristic surface for the equation (\ref{eq:2.1}),  i.e.
\begin{equation}                                     \label{eq:5.1}
\sum_{j,k=1}^ng^{jk}(x)S_{x_j}(x)S_{x_k}(x)=0\ \ \ \mbox{when}\ \ S(x)=0.
\end{equation}
Let $\Omega_{int}$  be the interior of $S$  and $\Omega_{ext}$  be the exterior of $S$.
The domain   $\Omega_{int}\times \R$ is called an artificial black hole if no signal emanating from
it  can reach $\Omega_{ext}\times\R$.  Analogously, $\Omega_{int}\times \R$  is an artificial white 
hole if no signal from $\Omega_{ext}\times\R$  can penetrate the interior of $S\times\R$.

Let $y$ be any point of $S$,  i.e. $S(y)=0$.

\begin{lemma}                                  \label{lma:5.1}
If $S\times\R$  is a characteristic surface  then 
$$
\sum_{j=1}^ng^{j0}(y)S_{x_j}(y)\neq 0.
$$
\end{lemma}

{\bf Proof:}
Since (\ref{eq:2.1})  is hyperbolic  the equation $\sum_{j.k=0}^n g^{jk}(y)\xi_j\xi_k=0$
has two distinct real roots $\xi_0^{(1)}(\xi),\xi_0^{(2)}(\xi)$ for any $\xi=(\xi_1,...,\xi_n)\neq 0$.
Taking $\xi=S_x(y)$ and using (\ref{eq:5.1}) we get 
$
g^{00}(y)\xi_0^2+2\sum_{j=1}^ng^{0j}(y)\xi_0S_{x_j}(y)=0.
$
Therefore $\xi_0^{(1)}=0,\xi_0^{(2)}=-2(g^{00}(y))^{-1}\sum_{j=1}^ng^{jk}(y)S_{x_j}(y)\neq 0$.
\qed

It follows from Lemma \ref{lma:5.1}  that either 
\begin{equation}                      \label{eq:5.2}
\sum_{j=1}^ng^{j0}(y)S_{x_j}(y)>0,\ \ S(y)=0,
\end{equation}
or
\begin{equation}                        \label{eq:5.3}
\sum_{j=1}^ng^{j0}(y)S_{x_j}(y)<0,\ \ S(y)=0.
\end{equation}
Denote by $K^+(y)\subset \R^{n+1}$  the half-cone
\begin{equation}                         \label{eq:5.4}
K^+(y)=\left\{(\xi_0,\xi_1,...,\xi_n):\sum_{j,k=0}^ng^{jk}(y)\xi_j\xi_k>0\right\}
\end{equation}
containing  $(1,0,...,0)$  and by $K^+(y)$ the dual half-cone 
\begin{equation}                         \label{eq:5.5}
K_+(y)=\left\{(\dot x_0,\dot x_1,...,\dot x_n)\in \R^{n+1}:\sum_{j,k=0}^ng_{jk}(y)\dot x_j\dot x_k>0,
\dot x_0>0\right\}.
\end{equation}
Since $K^+(y)$ and $K_+(y)$  are dual we have 
\begin{equation}                          \label{eq:5.6}
\sum_{j=0}^n\dot x_j \xi_j>0
\end{equation}
for any $(\dot x_0,...,\dot x_n)\in K_+(y)$  and any $(\xi_0,...,\xi_n)\in K^+(y)$.
We choose $S_x(y)$  to be the outward normal to $S$.
Assuming (\ref{eq:5.2}) we have $(\e,S_x(y))\in K^+(y)$  for any $\e>0$.
Using (\ref{eq:5.6}) and taking the limit when $\e\rightarrow 0$  we get that
$\sum_{j=1}^n\dot x_jS_{x_j}(y)\geq 0$  for all $(\dot x_0,...,\dot x_n)\in K_+(y)$,
i.e.
$\overline{K}_+(y)$ is contained in the half-space $\overline{P}_+(y)=
\left\{(\alpha_0,\alpha_1,...,\alpha_n):\sum_{j=1}^n\alpha_jS_{x_j}(y)\geq 0\right\}$.
In particular,  $K_+(y)$  is contained in the open half-space $P_+$.

A ray $x_0=x_0(s), x=x(s),s\geq 0,$
is called a forward time-like if $\left(\frac{dx_0(s)}{ds},\frac{dx(s)}{ds}  \right)\in
K_+(x(s))$  for all $s$.
It is known (c.f. [CH])  that the domain of influence of a point $(y_0,y)$  is the closure
of all forward time-like rays starting at $(y_0,y)$.  Therefore since $K_+(y)$  is
contained in the open half-space $P_+(y)$  for all $(y_0,y)\in S\times \R$  we have that the domain of
influence of $\Omega_{ext}\times \R$ is contained in $\overline{\Omega}_{ext}\times\R$,  i.e.
$\Omega_{int}\times\R$  is a white hole,  since no signal from $\Omega_{ext}\times\R$  may reach 
$\Omega_{int}\times\R$.

Consider now the case when (\ref{eq:5.3}) holds.  Then
$(\e,-S_x(y))\in K^+(y)$  for any $\e>0,y\in S$.  Therefore  passing to the limit 
when $\e\rightarrow 0$
we get that $K_+(y)$  is contained in the half-space ${P}_-(y)=
\left\{(\alpha_0,\alpha_1,...,\alpha_n):\sum_{j=1}^n\alpha_jS_{x_j}(y)< 0\right\}$.
Since  $S_x(y)$  is the outward normal to $S, y\in S$,  we get that the domain of influence of
$\Omega_{int}\times\R$  is contained in $\overline{\Omega}_{int}\times\R$, i.e.
$\Omega_{int}\times \R$ is a black hole.  We proved the following theorem:

\begin{theorem}                                     \label{theo:5.1}
Let $S\times\R$ be a characteristic surface for (\ref{eq:2.1}).  Then $\Omega_{int}\times\R$
is a white hole if (\ref{eq:5.2})
holds and a black hole if (\ref{eq:5.3}) holds.
\end{theorem}
In \S 8 and \S 9 we will give another proof of this theorem.

Let $\Delta(x)=\det[g^{jk}(x)]_{j,k=1}^n$.  
We assume that the surface $S_\Delta=\{x:\Delta(x)=0\}$  is a smooth closed surface.  Le $\Omega_{int}$
be the interior of $S_\Delta$  and $\Omega_{ext}$  be the exterior of $S_\Delta$.  We assume that
$\Delta(x)>0$  in $\overline{\Omega}\cap\Omega_{ext}$  and $\Delta(x)<0$  in 
$\overline{\Omega}\cap\Omega_{int}$.
Borrowing the terminology from  the general relativity we shall call $S_\Delta$ the ergosphere.
If $S_\Delta\times\R$  is a characteristic surface for (\ref{eq:2.1})  then 
$\Omega_{int}\times\R$  is a black hole if (\ref{eq:5.3}) holds and a white hole if (\ref{eq:5.2})
holds.   In the case of the Gordon equation the ergosphere has the form
$$
S_\Delta=\left\{x:\ |w(x)|^2=\frac{c^2}{n^2}\right\},
$$
and 
$$
g^{0j}(x)=\frac{(n^2(x)-1)cw_j(x)}{c^2-|w|^2}.
$$
If $S_\Delta\times \R$  is a characteristic surface then the normal to $S_\Delta$ is
colinear to $w(y)$  and $\Omega_{int}\times\R$  will be a black hole if  $w(y)$  is pointed 
inside  $\Omega_{int}$,  and $\Omega_{int}\times\R$ will be a white hole if $w(y)$  is pointed   
inside $\Omega_{ext}$.

Note that the black or white holes with the boundary $S_\Delta\times\R$  are not stable:  If we 
perturb slightly the metric $[g_{jk}(x)]_{j,k=0}^n$  then the ergoosphere changes slightly.  However
it will not  necessary    remain a characteristic surface and the black or the white hole will disappear.
In the next section we will find stable black and white holes.

\section{Stable black  and white holes.}
\label{section 6}
\init

Consider the case $n=2$,  i.e.  the case of two space variables $x=(x_1,x_2)$.
Let $S_\Delta$  be the ergosphere,  i.e. $\Delta(x)=g^{11}(x)g^{22}(x)-(g^{12}(x))^2=0$
on $S_\Delta$.  Suppose  $S_\Delta$ is a closed smooth curve and let  $S_1$  be
a smooth closed curve inside $S_\Delta$.  Denote by $\Omega_e$  the domain between $S_\Delta$
and $S_1$  and assume that $\Omega_e\subset \Omega$.  We shall call $\Omega_e$  the ergoregion.
We assume that $\Delta(x)<0$  on $\overline{\Omega}_e\setminus S_\Delta$  and that 
$S_\Delta$  is not characteristic at any $y\in S_\Delta$,  i.e.
\begin{equation}                                     \label{eq:6.1}
\sum_{j,k=1}^2g^{jk}(y)\nu_j(y)\nu_k(y)\neq 0,\ \ \forall y\in S_\Delta,
\end{equation}
where  $(\nu_1(y),\nu_2(y))$  is the nomal to $S_\Delta$.
Since $\Delta(x)<0$  in $\Omega_e$  we can define (locally)  two families of
characteristic curves $S^\pm(x)=\mbox{const}\ $
satisfying
\begin{equation}                                    \label{eq:6.2}
\sum_{j,k=1}^2g^{jk}(x)S_{x_j}^\pm(x)S_{x_k}^\pm(x)=0,\ \ x\in\Omega_e.
\end{equation}
It is shown in [E5] that there are 
two families $f^\pm(x)$ of vector fields such that $f^\pm(x)\neq (0,0)$  for $\forall x\in 
\overline{\Omega}_e, \ f^+(x)\neq f^-(x)$  for 
$x\in\overline{\Omega}_e\setminus S_\Delta,\ f^+(y)=f^-(y)$
on $S_\Delta$,  and $f^\pm(x)$  are tangent to $S^\pm(x)=\mbox{const}$.

Consider two systems of differential equations:
\begin{equation}                                          \label{eq:6.3}
\frac{d\hat x^+(\sigma)}{d\sigma}=f^+(\hat x^+(\sigma)),\ \ \sigma\geq 0,\ \ 
\hat x^+(0)=y\in S_\Delta,
\end{equation}  
\begin{equation}                                          \label{eq:6.4}
\frac{d\hat x^-(\sigma)}{d\sigma}=f^-(\hat x^-(\sigma)),\ \ \sigma\geq 0,\ \ 
\hat x^+(0)=y\in S_\Delta,
\end{equation}  
Note that $x=\hat x^\pm(\sigma)=y,\sigma\geq 0$,  are parametric equations
of characteristics (\ref{eq:6.2}).  It follows from (\ref{eq:6.1})  that
$f^+(y)=f^-(y)$  is  not tangent to $S_\Delta$  for all $y\in S_\Delta$.
Since the rank of $[g^{jk}(y)]_{j,k=1}^2$  on $S_\Delta$  is 1,  one can choose a smooth
vector $b(y),y\in S_\Delta$, such that
\begin{equation}                                            \label{eq:6.5}
\sum_{k=1}^2g^{jk}(y)b_k(y)=0,\ \ j=1,2.
\end{equation}
Note that $f^\pm(y)\cdot b(y)=0$.  We choose $f^\pm(y)$  to be pointed inside $S_\Delta$.

Consider the equations for the null-bicharacteristics:
\begin{equation}                                     \label{eq:6.6}
\frac{dx_j(s)}{ds}=2\sum_{k=0}^2g^{jk}(x(s))\xi_k(s),\ \ x_j(0)=y_j,\ \ 0\leq j\leq 2,
\end{equation}
\begin{equation}                                       \label{eq:6.7}
\frac{d\xi_p(s)}{ds}=-\sum_{j,k=0}^2g_{x_p}^{jk}((s))\xi_j(s)\xi_k(s),\ \xi_p(0)=\eta_p,\ \ 
0\leq p\leq 2.
\end{equation}
Here $x(s)=(x_1(s),x_2(s))$.  Since $g^{jk}(x)$  are independent of $x_0$  we have that 
$\xi_0(s)=\eta_s,\forall s$,  and
we choose $\eta_0=0$.

The bicharacteristic (\ref{eq:6.6}), (\ref{eq:6.7})  is a null-bicharacteristics if
$$
\sum_{j,k=0}^2g^{jk}(y)\eta_j\eta_k=0.
$$
Choosing $\eta_j=\pm b(y),\ 1\leq j\leq 2,\ \eta_0=0$  we get two null-bicharacteristics 
$x_0=x_0^\pm(s),  x=x^\pm(s),\ \xi_0=0,\xi=\xi^\pm(s)$
such that the projection of these null-bicharacteristics  on the  $(x_1,x_2)$-plane coincide 
with solutions $x=\hat x^\pm(\sigma)$  of the systems (\ref{eq:6.3}), (\ref{eq:6.4}),
i.e  $x=\hat x^\pm(\sigma),\sigma\geq 0$   and $x=x^\pm(s)$  are equal after a reparametrization 
$\sigma=\sigma^\pm(s),\ \frac{d\sigma^\pm(s)}{ds}>0$  for $s>0$.  We consider forward
null-bicharacteristics,  i.e. $\frac{dx_0^\pm(s)}{ds}>0$  for all $s$.  Therefore one can
take the time variable $x_0$  as a parameter on $x=x^\pm(\sigma)$.

The key observation in [E5]  is that for one of $x=\hat x^\pm(\sigma)$,  say for 
$x=\hat x^+(\sigma),\ \ \sigma=\sigma^+(s^+(x_0))$
increases  when $x_0$  increases,  and for $x=\hat x^-(\sigma),\ \ \sigma=\sigma^-(s^-(x_0))$
decreases  when $x_0$  increases.

Now we impose conditions on $S_1$  that will garantee the existence  of black and white
holes in $\Omega_e$.
We assume that $S_1$  is not characteristic.

Let $N(y)$  be the outward unit normal to $S_1, y\in S_1$.
Suppose that either 
\begin{description}
\item[$(a)$]
$\overline{K}_+(y)$  is contained in the open half-space $Q_+=
\{(\alpha_0,\alpha_1,\alpha_2):(\alpha_0,\alpha_1,\alpha_2)\cdot(0,N(y))>0,$
\\
or
\item[$(b)$]
$\overline{K}_+(y)$  is contained in the open half-space $Q_-=
\{(\alpha_0,\alpha_1,\alpha_2):(\alpha_0,\alpha_1,\alpha_2)\cdot(0,N(y))<0,$
\end{description}

{\bf Remark 6.1}
There  are equivalent forms of conditions $(a)$  and $(b)$.  Since 
 $\sum_{j,k=0}^ng_{jk}(x(s))\frac{dx_j(s)}{ds}\frac{dx_k(s)}{ds}=0$ 
for the null-bicharacteristics
 we have that
$\left(\frac{dx_0(s)}{ds}, \frac{dx_1(s)}{ds},\frac{dx_2(s)}{ds}\right)\in \overline{K}_+(y)$
for the forward null-bicharacteristic when $x(s_1)=y\in S_1$.  Therefore the condition $(a)$ is
equivalent to the condition:
\begin{description}
\item[($a_1$)]
The projection on $(x_1,x_2)$-plane  of all forward null-characteristics passing through
$y\in S_1$  leave $\Omega_e$  when $x_0$  increases.
\end{description}
Further,  the condition $(a_1)$  is equivalent to the following more simple condition:

Let $x=x^\pm(s)$  be the projection on $(x_1,x_2)$-plane of two forward null-bicharacteristics
such that $x=x^\pm(s)$ are the parametric equations of the characteristics $S^\pm(x)=\mbox{const}$,
i.e. $x=x^\pm(s)$  are solutions of  the differential equations (\ref{eq:6.3}), (\ref{eq:6.4})
after a reparametrization.  Assume that
\begin{description}
\item[($a_2$)]
$\frac{dx^\pm(s_1)}{ds}\cdot N(y)>0$  
when $x^\pm(s_1)=y$.
\end{description}
The condition $(a_2)$  follows from $(a_1)$.   The inverse is also true since the set 
of directions of the projections of all forward null-bicharacteristics passing through
$y$  is bounded by $\frac{dx^+(s_1)}{ds}$  and $\frac{dx^-(s_1)}{ds}$.

Conditions $(b_1),\ (b_2)$  are similar to $(a_1),(a_2)$ 
when the sign of the inner product in  $(a)$  is negative.

\begin{theorem}(c.f. [E5])                                      \label{theo:6.1}
Let $\partial\Omega_e=S_\Delta\cup S_1$,  where $S_\Delta$  is the ergosphere,  i.e. $\Delta(y)=0$  
on $S_\Delta$.
Suppose (\ref{eq:6.1})  holds on $S_\Delta$ and either (a) or (b)  hold on $S_1$.  Then 
there exists a closed Jordan curve $S_0(x)=0$   inside $\Delta_e$  such that $S_0\times \R$  is
the boundary of either black or white hole.
\end{theorem}

The proof of Theorem \ref{theo:6.1}  is based on the Poincare-Bendixson theorem  (c.f. [H]).  
Suppose $(a)$  holds.  Then 
the solution of (\ref{eq:6.4}) cannot reach  $S_1$.  Indeed,  suppose $\hat x^-(\sigma_1)=y_1\in S_1$
for some $\sigma_1>0$.  Then $\hat x^-(\sigma)$  leaves $\Omega_e$  when $\sigma>\sigma_1$.
From other side, when $\sigma$ increases $x_0$ decreases.  Therefore  $x=\hat x^-(\sigma^-(x_0))$  leaves 
$\Omega_e$  when $x_0$  decreases,  and this contradicts the condition $(a_1)$.  Since $x=\hat x^-(\sigma)$
never reaches $S_1$  the limit set of the trajectory $x=\hat x^-(\sigma)$  is contained inside $\Omega_e$.
Then by the Poincare-Bendixson theorem there exists a limit cycle $S_0(x)=0$,
i.e. a Jordan curve that is a periodic solution of $\frac{d\hat x}{d\sigma}=f^-(\hat x(\sigma))$.  
Therefore 
$S_0\times \R$  is a characteristic surface and it is the boundary of a black or a white hole.
In the case when the condition $(b)$  holds we have that the solution of (\ref{eq:6.3})
never reach $S_1$.   Therefore again by the Poincare-Bendixson Theorem there exists a  
black or white hole.

Applying Theorem \ref{theo:6.1} to the Gordon equation we get
\begin{theorem}                                \label{theo:6.2}
Let $S_\Delta$  be the ergosphere,  i.e. $|w|^2=\frac{c^2}{n^2(x)}$.  Suppose $w(x)$  is not 
colinear with the normal to $S_\Delta$  for any $x\in S_\Delta$.   Suppose that either
\begin{equation}                               \label{eq:6.8}
(n^2(x)-1)^{\frac{1}{2}}(v(x)\cdot N(x))>1\ \ \mbox{on}\ \ S_1,
\end{equation}
or
\begin{equation}                               \label{eq:6.9}
(n^2(x)-1)^{\frac{1}{2}}(v(x)\cdot N(x))<-1\ \ \mbox{on}\ \ S_1,
\end{equation}
where $v(x)=\left(1-\frac{|w|^2}{c^2}\right)^{-\frac{1}{2}}\frac{w(x)}{c},\ 
N(x)$  is the outward unit normal to $S_1$.

Then there exists a limit cycle $S_0(x)=0$   and $S_0\times\R$  is the boundary of a black 
or a white hole.
\end{theorem}

{\bf Remark 6.1}
Note that the black or white holes obtained by Theorems \ref{theo:6.1} and \ref{theo:6.2}
are stable since the assumptions remain valid when we slightly deform the metric.

\section{Rotating black holes. Examples.}
\label{section 7}
\init

{\bf Example 1 ([V]).  Acoustic black hole.}
Consider a fluid flow with velocity field
\begin{equation}                            \label{eq:7.1}
v=(v^1,v^2)=\frac{A}{r}\hat r +\frac{B}{r}\hat\theta,
\end{equation}
where $r=|x|,\hat r=\left(\frac{x_1}{|x|},\frac{x_2}{|x|}\right),        
\ \hat \theta=\left(-\frac{x_2}{|x|},\frac{x_1}{|x|}\right),\ A $  and $B$ are constants.
The inverse of the metric tensor has the following form in this case:
\begin{eqnarray}                                     \label{eq:7.2}
g^{00}=\frac{1}{\rho c},\ g^{oj}=g^{j0}=\frac{1}{\rho c}v^j,\ \ 1\leq j\leq 2,
\\
\nonumber
g^{jk}=\frac{1}{\rho c}(-c^2\delta_{ij}+v^jv^k),\ \ 1\leq j,\ k\leq 2,
\end{eqnarray}
where $c$  is the sound speed,   $\rho$  is the density.

Consider the case $A>0,\ B>0$.  Assume $\rho=c=1$.  Then 
the ergosphere is $r=\sqrt{A^2+b^2}$.  Consider the domain
$\Omega_e=\{r_1\leq r\leq \sqrt{A^2+B^2}\}$,  where $r_1<A$.  In polar coordinates 
$(r,\theta)$  the differential equations (\ref{eq:6.3}), (\ref{eq:6.4}) have the form:
\begin{equation}                                   \label{eq:7.3}
\frac{dr}{ds}=A^2-r^2,\ \ \frac{d\theta}{ds}=\frac{AB}{r}+\sqrt{A^2+B^2-r^2},
\end{equation}
and
\begin{equation}                                   \label{eq:7.4}
\frac{dr}{ds}=-1,\ \ \frac{d\theta}{ds}=\frac{1-\frac{B^2}{r^2}}{\frac{AB}{r}+\sqrt{A^2+B^2-r^2}}.
\end{equation}
We have that $r=A$  is a limit cycle and $\{r=A\}\times\R$  is the boundary of a white hole (c.f. [E5]).

{\bf Example 2}
Consider a fluid flow with the velocity 
$$
v=A(r)\hat r +B(r)\hat\theta,
$$
where $r_1\leq r\leq r_0,\ A(r),B(r)$  are smooth, $B(r)>0,\ A^2(r_0)+B^2(r_0)=1,\ 
A^2(r)+B^2(r)>1$  on $[r_1,r_0),\ A(r)+1$  has simple zeros $\alpha_1,...,\alpha_{m_1}$   on
$(r_1,r_0),\ A(r)-1$  has simple zeros $\beta_1,...,\beta_{m_2}$   on
$(r_1,r_0),\ \alpha_j\neq\beta_k,\ \forall j,\ \forall k,\ |A(r_1)|>1$.
Here $r=r_0$  is the ergosphere.  The differential equations (\ref{eq:6.3}), (\ref{eq:6.4})
have the following form in polar coordinates $(r,\theta)$:
\begin{equation}                                       \label{eq:7.5}
\frac{dr}{ds}=A(r)-1,\ \ \ \frac{d\theta}{ds}=\frac{A(r)B(r)+\sqrt{A^2(r)+B^2(r)-1}}{A(r)+1},
\end{equation}
and
\begin{equation}                                       \label{eq:7.6}
\frac{dr}{ds}=A(r)+1,\ \ \ \frac{d\theta}{ds}=\frac{A(r)B(r)-\sqrt{A^2(r)+B^2(r)-1}}{A(r)-1}.
\end{equation}
Here $r=\alpha_j,\ 1\leq j\leq m_1$,  and $r=\beta_n,\ 1\leq k\leq m_2$  are limit cycles and there 
are $m_1+m_2$  black and white holes.

{\bf Axially symmetric metrics.}

Consider the equation (\ref{eq:2.1}) in $\Omega\times\R$  where $\Omega$  is a three-dimentional domain.
Let $(r,\theta,\varphi)$  be the spherical coordinates in $\R^3$.  Suppose 
$g^{jk}$  are independent of $\varphi$.

Consider a characteristic surface $S$  independent of $\varphi$  and $x_0$,  i.e. $S$
depends on $r$ and $\theta$  only.   Then $S$  satisfies an equation
\begin{equation}                                              \label{eq:7.7}
a^{11}(r,\theta)\left(\frac{\partial S}{\partial r}\right)^2 
+ 2 a^{12}(r,\theta)\frac{\partial S}{\partial r}\frac{\partial S}{\partial \theta}+
a^{22}(r,\theta)\left(\frac{\partial S}{\partial \theta}\right)^2=0.
\end{equation}
We assume that $a^{ij}(r,\theta)$  are also independent of $\varphi,\ 1\leq j,k\leq 2.$

Consider (\ref{eq:7.7}) in two-dimensional domain $\omega$  where 
$\delta_1\leq r\leq \delta_2,\ 0<\delta_3<\theta<\pi-\delta_4$   when
$(r,\theta)\in\omega$.
Here $\delta_j>0,\ 1\leq j\leq 4.$

Assuming that $\omega$  and $a^{jk}(r,\theta),\ 1\leq j,k\leq 2,$
satisfy the condition of Theorem \ref{theo:6.1},  we can prove the existence of black or white holes
whose boundary is $S_0\times S^1\times\R$,   where $\varphi\in S^1,\ x_0\in \R$  and $S_0$
is a Jordan curve in $\omega$.

Such black (white) holes are called the rotating black (white) holes.

\section{Black holes and inverse problems I.}
\label{section 8}
\init

In this section we consider the case of the black or the white hole bounded by $S_\Delta\times \R$,
where $S_\Delta$  is  the ergosphere.  Suppose $\Omega_{int}$  is a black hole,  i.e. 
(\ref{eq:5.3})  holds.
Let $L$  be the operator (\ref{eq:2.1}).  Consider $u(x_0,x)$  in $(\Omega\cap\Omega_{ext})\times(0,T)$  
such that
\begin{equation}                                   \label{eq:8.1}
Lu=f,\ \ \ (x_0,x)\in(\Omega\cap\Omega_{ext})\times(0,T),
\end{equation}
\begin{equation}                                   \label{eq:8.2}
u(0,x)=0,\ \ \ \frac{\partial u(0,x)}{\partial x_0}=0,\ \ \ x\in (\Omega\cap\Omega_{ext}),
\end{equation}
\begin{equation}                                  \label{eq:8.3}
u{\Huge|}_{\partial\Omega_0\times(0,T)}=g.
\end{equation}
We do not impose  
any boundary 
conditions on $S_\Delta\times(0,T)$  and we assume, for the simplicity, that there is no
obstacles between $\partial\Omega_0$  and $S_\Delta$.
We shall prove an estimate of $u(x_0,x)$  in terms of $g$ and $f$.

Denote $Hu=\sum_{j=1}^ng^{0j}(x)u_{x_j}$.  Consider the equality
$$
(Lu,g^{00}u_{x_0}+Hu)=(f,g^{00}u_{x_0}+Hu),
$$
where $(u,v)$  is the inner product in $L_2((\Omega\cap\Omega_{int})\times(0,T))$.
We shall denote by $Q_p(u,v))$  for $p\geq 1$  the expression of the form:
\begin{equation}                                     \label{eq:8.4}
(Q_pu,v)
=\int_0^T\int_{\Omega\cap\Omega_{int}}\sum_{j,k=0}^nq_{jkp}(x)u_{x_j}v_{x_k}dxdx_0.
\end{equation}
Denote
\begin{eqnarray}                                  \label{eq:8.5}
I_1=
\left(
\sum_{j=1}^n\frac{1}{\sqrt{|g|}}
\frac{\partial}{\partial x_j}\left(\sqrt{|g|}g^{0j}\frac{\partial u}{\partial x_0}\right)
\right.
\nonumber
\left.+\sum_{j=0}^n\frac{1}{\sqrt{|g|}}
\frac{\partial}{\partial x_0}\left(\sqrt{|g|}g^{0j}\frac{\partial u}{\partial x_j}\right),
g^{00}u_{x_0}+Hu\right)
\\
\ \ \ \ \ \ \ \ \ \ \ \ \ \ \ \ \ 
\stackrel{def}{=}I_{11}+I_{12}+I_{13} \ \ \ \ \ \ \ \ \ \ \ \ \ \ \ \ \ \ \ \ \ \ \ \ \ \ \ \
\end{eqnarray}
We have
\begin{eqnarray}                                     \label{eq:8.6} 
I_{11}=\left(\sum_{j=0}^n\frac{1}{\sqrt{|g|}}
\frac{\partial}{\partial x_0}\left(\sqrt{|g|}g^{0j}\frac{\partial u}{\partial x_j}\right),
g^{00}u_{x_0}+Hu\right)
\\
\nonumber
=\frac{1}{2}\int_{\Omega\cap\Omega_{int}}\left(\sum_{j=0}^ng^{0j}u_{x_j}\right)^2dx\left|_0^T\right.+
Q_1(u,u),
\end{eqnarray}
where $a|_0^T$  means $a(T)-a(0)$.  Note that
\begin{eqnarray}                                        \label{eq:8.7}
I_{12}=
\left(\sum_{j=1}^n\frac{1}{\sqrt{|g|}}
\frac{\partial}{\partial x_j}\left(\sqrt{|g|}g^{0j}\frac{\partial u}{\partial x_0}\right),
Hu\right)
\\
\nonumber
=\frac{1}{2}\int_{\Omega\cup\Omega_{ext}}(Hu)^2dx\left|_0^T\right.+Q_2(u,u).
\end{eqnarray}
Also we have
\begin{eqnarray}                                        \label{eq:8.8}
I_{13}=
\left(\frac{1}{\sqrt{|g|}}
\sum_{j=1}^n
\frac{\partial}{\partial x_j}\left(\sqrt{|g|}g^{0j}\frac{\partial u}{\partial x_0}\right),
g^{00}u_{x_0}\right)
\\
\nonumber
=\frac{1}{2}\int_{\Omega\cup\Omega_{ext}}\sum_{j=1}^n\frac{\partial}{\partial x_j}(g^{0j}g^{00}u_{x_0}^2)
dxdx_0+Q_3(u,u).
\end{eqnarray}
By the divergence theorem  we get
\begin{eqnarray}                                          \label{eq:8.9}
I_{13}=-\int_0^T\int_{S_\Delta}
\frac{1}{2}\left(\sum_{j=1}^ng^{0j}(x)\nu_j(x)\right) g^{00}u_{x_0}^2dsdx_0
\\
\nonumber
+\int_0^T\int_{\partial\Omega_0}
\frac{1}{2}\left(\sum_{j1}^ng^{0j}(x)N_j(x)\right)g^{00}u_{x_0}^2dsdx_0+Q_3(u,u),
\end{eqnarray}
where 
$ds$  is the area element on $S_\Delta$  and $\partial\Omega_0$,  respectively,
$N(x)=(N_1,...,N_n)$  is the outward unit normal to $\partial\Omega_0$ 
and $\nu=(\nu_1,...,\nu_n)$  is an outward normal to $S_\Delta$.  
Note that $\nu$  is the inward normal with respect to $\Omega\cap\Omega_{ext}$
Therefore
\begin{eqnarray}                                       \label{eq:8.10}
I_1=I_{11}+I_{12}+I_{13}\ \ \ \ \ \ \ \ \ \ \ \ \ \ \ \ \ \ \ \ \ \ \ \ \ \ \ \ \ \ \ \
\\
\nonumber
=\frac{1}{2}\int_{\Omega\cap\Omega_{ext}}
[(\sum_{j=0}^ng^{j}u_{x_j}
)^2+(Hu)^2]dx\left|_0^T\right.
\\
\nonumber
-\int_0^T\int_{S_\Delta}
\frac{1}{2}\left(\sum_{j=1}^ng^{0j}(x)\nu_j(x)\right) g^{00}u_{x_0}^2dsdx_0
\\
\nonumber
+\int_0^T\int_{\partial\Omega_0}
\frac{1}{2}\left(\sum_{j=1}^ng^{0j}(x)N_j(x)\right)g^{00}u_{x_0}^2dsdx_0+Q_4(u,u).
\end{eqnarray}

Now consider
$$
I_2=\left(\sum_{j,k=1}^n\frac{1}{\sqrt{|g|}}
\frac{\partial}{\partial x_j}\left(\sqrt{|g|}g^{jk}\frac{\partial u}{\partial x_k}\right),
g^{00}u_{x_0}\right).
$$
Integrating  by parts in $x_j$  and taking tnto account that
$$
-\sum g^{jk}u_{x_k}u_{x_0x_j}=
-\frac{1}{2}
\frac{\partial}{\partial x_0}
\left(\sum_{j,k=1}^n g^{jk}u_{x_j}u_{x_k}\right),
$$
we get:
\begin{eqnarray}                                           \label{eq:8.11}
I_2=-\frac{1}{2}\int_{\Omega\cap\Omega_{ext}}g^{00}
\sum_{j,k=1}^ng^{jk}(x)u_{x_j}u_{x_k}dx\left|_0^T\right.
\\
\nonumber
-\int_0^T\int_{S_\Delta}\left(\sum_{j,k=1}^n g^{jk}u_{x_j}(x)\nu_k(x)g^{00}u_{x_0}\right)dsdx_0
\\
\nonumber
+\int_0^T\int_{\partial\Omega_0}
\left(\sum_{j,k=1}^ng^{jk}(x)u_{x_j}N_k(x)g^{00}u_{x_0}\right)dsdx_0
+Q_5(u,u).
\end{eqnarray}
Since $S_\Delta$  is an ergosphere and a characteristic surface we have (c.f. (\ref{eq:6.5})):
\begin{equation}                                            \label{eq:8.12}
\sum_{k=1}^ng^{jk}(x)\nu_k(x)=0\ \ \ \mbox{on}\ \ S_\Delta,\ j=1,...,n.
\end{equation}                            
Let
\begin{equation}                                      \label{eq:8.13}
I_3=\left(\sum_{j,k=1}^n
\frac{\partial}{\partial x_j}\left(\sqrt{|g|}g^{jk}\frac{\partial u}{\partial x_k}\right),
Hu\right).
\end{equation}
Integraiting by parts in $x_j,\ 1\leq j\leq n,$  we get
\begin{eqnarray}                                     \label{eq:8.14}
\ \ \ \ \ \ \ \ I_3=-\left(\sum_{j,k=1}^ng^{jk}u_{x_k},Hu_{x_j}\right)
-\int_0^T\int_{S_\Delta}\sum_{j,k=1}^ng^{jk}\nu_ju_{x_k}Hu\ dsdx_0
\\
\nonumber
+\int_0^T\int_{\partial \Omega_0}\sum_{j,k=1}^ng^{jk}N_ju_{x_k}Hu\ dsdx_0+Q_6(u,u).
\end{eqnarray}
We have 
\begin{eqnarray}                                       \label{eq:8.15}
\int_0^T\int_{\Omega\cap\Omega_{ext}}
\frac{\partial}{\partial x_p}g^{0p}\sum_{j,k=1}^ng^{jk}u_{x_j}u_{x_k} dxdx_0\ \ \ \ \ \ \ \ \ \ \ \
\ \ \ \ \ \  
\\
\nonumber
=-\int_0^T\int_{S_\Delta}g^{0p}\nu_p\sum_{j,k=1}^ng^{jk}u_{x_j}u_{x_k}dsdx_0
+
\int_0^T\int_{\partial\Omega_0}g^{0p}N_p\sum_{j,k=1}^ng^{jk}u_{x_j}u_{x_k}dsdx_0
\end{eqnarray}
Using (\ref{eq:8.15}) we obtain
\begin{eqnarray}                                       \label{eq:8.16}
\nonumber
-\left(\sum_{j,k=1}^ng^{jk}u_{x_k},Hu_{x_j}\right)=
\int_0^T\int_{S_\Delta}\left(\sum_{p=1}^ng^{0p}\nu_p\right)\sum_{j,k=1}^ng^{jk}u_{x_j}u_{x_k}dsdx_0
\\
-\int_0^T\int_{\partial\Omega_0}\left(\sum_{p=1}^ng^{0p}N_p\right)
\sum_{j,k=1}^ng^{jk}u_{x_j}u_{x_k}dsdx_0
+Q_7(u,u).
\end{eqnarray}
Note that the first integral in (\ref{eq:8.16})  is nonnegative since 
$\sum_{p=1}^n g^{0p}\nu_p <0$ on $S_\Delta$  and 
$\sum_{j,k=1}^ng^{jk}u_{x_j}u_{x_k}\leq 0$  in $\Omega\cap \overline\Omega_{int}$  since 
the matrix $[g^{jk}]_{j,k=1}^n$  has one zero eigenvalue and $n-1$  negative eigenvalues
on $S_\Delta$.

Now we shall estimate the integrals over  $\partial\Omega_0\times (0,T)$.
Let $\alpha(x)\in C_0^\infty(\R^n),\ \alpha(x)=1$  near $\partial\Omega_0$  and 
$\alpha(x)=0$  near  $S_\Delta$.   We have from (\ref{eq:8.1}), (\ref{eq:8.2}), (\ref{eq:8.3}) 
that $v=\alpha u$  satisfies
\begin{equation}                                      \label{eq:8.17}
Lv=\alpha f+L_1u \ \ \ \mbox{in}\ \ \Omega_0\times (0,T),
\end{equation}
\begin{equation}                                      \label{eq:8.18}
v|_{x_0=0}=0,\ \ v_{x_0}|_{x_0=0}=0,\ \ x\in \Omega_0,
\end{equation}
\begin{equation}                                      \label{eq:8.19}
v|_{\partial\Omega_0\times (0,T)}=g.
\end{equation}
Here $g$  is the same as in (\ref{eq:8.3})  and $\mbox{ord}\ L_1\leq 1$.

Since $L$ is strictly hyperbolic and $\partial\Omega_0$  ia not characteristic the following 
estimate for the solution of 
(\ref{eq:8.17}), (\ref{eq:8.18}), (\ref{eq:8.19})
holds (c.f.  for example,  [Ho],  see also [E7]):
\begin{eqnarray}                                  \label{eq:8.20}
\|v_{x_0}(T,\cdot)\|_0^2+\|v(T,\cdot)\|_1^2+\left[\frac{\partial v}{\partial N}\right]_0^2
\ \ \ \ \ \ \ \ \ \ \ \ \ \ 
\\
\nonumber
\leq C_T \left([g]_1^2
+\int_0^T\|f(x_0,\cdot)\|_0^2dx_0 
+\int_0^{T}(\|u(x_0,\cdot)\|_1^2+\|u_{x_0}(x_0,\cdot)\|_0^2)dx_0\right),
\end{eqnarray}
where $[w]_m$  is the norm in $H^m(\partial\Omega_0\times(0,T))$  and $\|w(x_0,\cdot)\|_m$  is the norm
in $H^m(\Omega_0),\ x_0$  is fixed.
All integrals over $\partial\Omega_0\times(0,T)$ in  $I_2$  and $I_3$  have the form
$$
I_4=\int_0^T \int_{\partial\Omega_0}\sum_{j,k=0}^np_{jk}u_{x_j}u_{x_k}dsdx_0.
$$
Therefore
\begin{eqnarray}                                       \label{eq:8.21}
|I_4|\leq C\left([g]_1^2+\left[\frac{\partial v}{\partial N}\right]_0^2
\right) \ \ \ \ \ \ \ \ \ \ \ \ \ \ \ \ \
\\
\nonumber
\leq C\left([g]_1^2+\int_0^T\|f(x_0,\cdot)\|_0^2dx_0+ \int_0^T(\|u(x_0,\cdot)\|_1^2+
\|u_{x_0}(x_0,\cdot)\|_0^2)dx_0\right),
\end{eqnarray}
where we used (\ref{eq:8.20})   to  estimate $\left[\frac{\partial v}{\partial N}\right]_0^2$.
Note that the norm 
\begin{equation}                                      \label{eq:8.22}
\int_{\Omega\cap\Omega_{ext}}[(u_{x_0}+Hu)^2+(Hu)^2-\sum_{j,k=1}^ng^{jk}u_{x_j}u_{x_k}]dx
\end{equation}
is equivalent to $\|u(x_0,\cdot)\|_1^2+\|u_{x_0}(x_0,\cdot)\|_0^2$.  Also
\begin{equation}                                         \label{eq:8.23}
|(f,g^{00}u_{x_0}+Hu)|\leq \frac{1}{2}\int_0^T\|f(y_0,\cdot)\|_0^2dy_0+
C\int_0^T (\|u_{x_0}\|_0^2+\|u(y_0,\cdot)\|_1^2)dy_0.
\end{equation}
Combining 
(\ref{eq:8.10}), (\ref{eq:8.11}), (\ref{eq:8.14}), (\ref{eq:8.16}), (\ref{eq:8.21}), taking into
account (\ref{eq:8.12})
and applying all estimates to the interval $(0,t_0)$  instead of $(0,T),\ t_0\leq T$,
we get
\begin{eqnarray}                                   \label{eq:8.24}
C(\|u(t_0,\cdot)\|_1^2                          
+\|u_{x_0}(t_0,\cdot)\|_0^2)
\ \ \ \ \ \ \ \ \ \ \ \ \ \ \ \ \ \ \ 
\ \ \ \ \ \ \ \ \ \ \ \ \ \ \ \ \ \ \ 
\\
-\int_0^T\int_{S_\Delta}\left(\sum_{j=1}^ng^{0j}\nu_j(x)\right)
\left(g^{00}u_{x_0}^2-\sum_{j,k=1}^ng^{jk}u_{x_j}u_{x_k}\right)dsdx_0  
\nonumber
\\
\leq C\int_0^T\|f(y_0,\cdot)\|_0^2dy_0+C\int_0^{t_0}(\|u(y_0,\cdot)\|_1^2
+\|u_{x_0}(y_0,\cdot)\|_0^2)dy_0+C[g]_1^2.
\nonumber
\end{eqnarray}
Note that the inequality $b(t)\leq C\int_0^t b(\tau)d\tau +d$  implies $b(t)\leq e^{ct}d$.  
Therefore we get
\begin{eqnarray}                                   \label{eq:8.25}
C(\|u(t_0,\cdot)\|_1^2                          
+\|u_{x_0}(t_0,\cdot)\|_0^2)
\ \ \ \ \ \ \ \ \ \ \ \ \ \ \ \ \ \ \ 
\ \ \ \ \ \ \ \ \ \ \ \ \ \ \ \ \ \ \ 
\\
-\int_0^T\int_{S_\Delta}\left(\sum_{j=1}^ng^{0j}\nu_j(x)\right)
\left(g^{00}u_{x_0}^2-\sum_{j,k=1}^ng^{jk}u_{x_j}u_{x_k}\right)dsdx_0  
\nonumber
\\
\leq C\int_0^T\|f(y_0,\cdot)\|_0^2dy_0
+C[g]_1^2,
\nonumber
\end{eqnarray}
where $0\leq t_0\leq T$.

We proved the following theorem:
\begin{theorem}                                       \label{theo:8.1}
Let $u(x_0,x)$  be the solution of 
(\ref{eq:8.1}), (\ref{eq:8.2}), (\ref{eq:8.3})
in $\Omega\cap \Omega_{ext}$.  
Let the ergosphere $S_\Delta$  be a characteristic surface and (\ref{eq:5.3})  holds.  Then 
$u(x_0,x)$  satisfies 
(\ref{eq:8.25}).
\end{theorem}
Note that $g^{00}u_{x_0}^2-\sum_{j,k=1}^nu_{x_j}u_{x_k}\geq 0$  
on $S_\Delta\times[0,T]$ 
since $S_\Delta$  is an ergosphere.

Theorem \ref{theo:8.1}  implies that the domain of dependence of 
$(\Omega\cap\Omega_{ext})\times\R$ 
is contained in $(\overline\Omega\cap\overline\Omega_{ext})\times\R$.
 Suppose $u$  is a solution of (\ref{eq:2.1}) and 
$\mbox{supp\ } u\subset\overline\Omega_{int}$ for $x_0\leq t_0$,  i.e. $u=0$  for 
$x\in(\Omega_{ext}\cap\Omega)\times(-\infty,t_0)$.  Then Theorem \ref{theo:8.1}
implies that $u=0$ for  $(\overline\Omega\cap\overline\Omega_{ext})\times[t_0,+\infty),$
i.e. $\mbox{supp\ }u\subset\overline\Omega_{int}\times\R$.
Therefore the domain of influence of $\Omega_{int}\times\R$  is contained in 
$\overline\Omega_{int})\times\R$,  i.e.  is a black hole.

Now we shall discuss the nonuniqueness of the inverse problem in the presence of a black hole.
Consider two initial-boundary value problem (\ref{eq:2.1}), (\ref{eq:2.5}), (\ref{eq:2.6})
for the operators $L_1$  and $L_2$  that differ only in $\Omega_{int}$.  Since $L_1=L_2$  in 
$\Omega_{ext}$  and we assume $f$  is the same for $L_1$  and  $L_2$  we get by Theorem \ref{theo:8.1}
that $u_1=u_2$  in $(\Omega\cap\Omega_{ext})\cap\R$  where $u_1$  and $u_2$  are the solutions
of the corresponding nitial-boundary value problems.
Therefore $\Lambda_1=\Lambda_2$  
on $\partial\Omega_0\times\R$.  Since $L_1\neq L_2$  in $\Omega_{int}$,  i.e.
we have a nonuniquenes of the inverse problem.
\qed

Consider now the case when  $S_\Delta$  is a characteristic surface and  (\ref{eq:5.2}) holds.
Suppose $\Omega_{int}\cap\Omega$  contains an obstacle $\Omega_1$ (it may be 
no obstacles or more than one obstacle,  but we consider the case of one obstacle for 
the definiteness).   Integrating by parts as in the proof of Theorem \ref{theo:8.1} and
using (\ref{eq:5.2})  instead of (\ref{eq:5.3}) we get the following theorem:
\begin{theorem}                                   \label{theo:8.2}
Consider the initial-boundary value problem:
\begin{equation}                    \label{eq:8.26}
Lu=f
\ \ \ \mbox{in}\ \ \ (\Omega_{int}\cap\Omega)\times (0,T),
\end{equation}
\begin{equation}                                  \label{eq:8.27}
u|_{\partial\Omega_1\times (0,T)}=g,
\ \ u(0,x)=0,\ \ u_{x_0}(0,x)=0\ \ \ \mbox{in}\ \ \Omega\cap\Omega_{int}.
\end{equation}
Suppose that the ergosphere  $S_\Delta$  is a characteristic surface and (\ref{eq:5.2}) holds.
Then an estimate of the form (\ref{eq:8.25})
holds in $(\Omega_{int}\cap\Omega)\times(0,T)$  with the following changes:  Integral over
$S_\Delta\times(0,T)$  must be taken with plus sign,  $\|u\|_s$  are the norms in 
$H^s(\Omega_{int}\cap\Omega),\ [g]_1$  is
the norm in $H^1(\partial\Omega_1\times(0,T))$.
\end{theorem}

The consequence of Theorem \ref{theo:8.2} is that the domain of dependence of 
$(\Omega_{int}\cap\Omega)\times\R$
is contained  in
$(\overline\Omega_{int}\cap\overline\Omega)\times\R$.  Therefore if $u(x_0,x)$
is the solution of (\ref{eq:2.1}) 
and $\mbox{supp\ }u\subset\overline\Omega_{ext}\cap\overline\Omega$
for $x_0\leq t_0$  then 
$\mbox{supp\ }u\subset\overline\Omega_{ext}\cap\overline\Omega$
for all $x_0 > t_0$,  i.e.  $\Omega_{int}\times\R$  is a white hole.
If $u(x_0,x)$  is the solution of  (\ref{eq:2.1}), (\ref{eq:2.5}), (\ref{eq:2.6}),
then $u(0,x)=0$  in $\Omega_{int}\cap\Omega,\ u|_{\partial\Omega_1\times\R}=0$.
Then by Theorem \ref{theo:8.2}  $u=0$  in $(\Omega_{int}\cap\Omega)\times\R$.
Therefore we can change the coefficients of $L$ in $\Omega_{int}\cap\Omega$
without changing the solution of
(\ref{eq:2.1}), (\ref{eq:2.5}), (\ref{eq:2.6}),
i.e. in the case of a white hole we again have a nonuniqueness of the solution of the
inverse problem.
\qed

Let $u(x_0,x)$  be the solution of
\begin{equation}                               \label{eq:8.28}
Lu=f\ \ \mbox{in}\ \ \Omega_{ext}\times(0,T),
\end{equation}
\begin{equation}                               \label{eq:8.29}
u(0,x)=\varphi_0(x),\ \ u_{x_0}(0,x)=\varphi_1(x),\ \ x\in\Omega_{ext},
\end{equation}
i.e. we consider 
(\ref{eq:8.28}), (\ref{eq:8.29})  in unbounded domain 
$\Omega_{ext}=\R^n\setminus\overline\Omega_{int}$.
We assume
that $[g^{jk}(x)]_{j,k=1}^n$  are smooth and have bounded derivatives of any order,
$g^{00}(x)\geq C_0>0,
\ \sum_{j,k=1}^ng^{jk}(x)\xi_j\xi_k\leq -C_0\sum_{j=1}^n\xi_j^2,\ x\in\overline\Omega_{ext}$
is  large,
and (\ref{eq:5.3}) holds.
Repeating the proof of Theorem \ref{theo:8.1} (with the simplification that we do not have 
the boundary condition (\ref{eq:8.3})),  we get  for any $T>0$:
\begin{eqnarray}                                \label{eq:8.30}
\max_{0\leq x_0\leq T}(\|u(x_0,\cdot)\|_{1,\Omega_{ext}}^2+
\|u_{x_0}(x_0,\cdot)\|_{0,\Omega_{ext}}^2)
\\
\nonumber
-\int_0^T\int_{S_\Delta}
(u_{x_0}^2-\sum_{j,k=1}^ng^{jk}u_{x_j}u_{x_k}
)
(\sum_{j=1}^ng^{j0}\nu_j
)dsdy_0
\\
\nonumber
\leq C([\varphi_0]_{1,\Omega_{ext}}^2+[\varphi_1]_{0,\Omega_{ext}}^2)
+C\int_0^T\|f(x_0,\cdot)\|_{0,\Omega_{ext}}^2dx_0.
\end{eqnarray}
Therefore the following theorem holds:
\begin{theorem}                                      \label{theo:8.3}
Suppose $u(x_0,x)$  satisfies 
(\ref{eq:8.28}), (\ref{eq:8.29}).  Suppose the ergosphere $S_\Delta$  is
a characteristic surface and (\ref{eq:5.3})  holds.
Then the estimate (\ref{eq:8.30})  holds for any $T>0$.
\end{theorem}

A consequence of Theorem \ref{theo:8.3}  is that 
$D_+(\Omega_{int}\times\R)\subset\overline\Omega_{int}\times\R$, i.e.
 $\Omega_{int}\times\R$  is a black hole.
\qed

An important problem is the determination of black or white holes by the boundary measurements on
$\partial\Omega_0\times(0,T)$  or on $\Gamma\times(0,T_0)$  where $\Gamma$  is an open part of
$\partial\Omega_0$.  Let $S_\Delta$  be the ergosphere inside $\Omega_0$.  
It does not matter in this subsection whether
 $S_\Delta$  forms  a black (or a white) hole.  The following theorem
is straightforward application of the proof of Theorem \ref{theo:2.1}:
\begin{theorem}                                   \label{theo:8.4}
The DN operetor $\Lambda$  given on $\Gamma\times(0,+\infty)$  determines $S_\Delta$  
up to a diffeomorphism (\ref{eq:2.8}).
\end{theorem} 
We note that the determination of $S_\Delta$  requires to take measurements on 
$\Gamma\times(0,+\infty)$.
It is not enough to know the Cauchy data on $\Gamma\times(0,T)$  for any finite $T$.  
The explanation of
this phenomenon is the following:
The proof of Theorem \ref{theo:2.1} allows to recover metric tensor $[g^{jk}]$  (up
to a diffeomorphism) gradually starting from the boundary $\partial\Omega_0$.  
The 
recovery of the metric at some point $x^{(1)}$  inside $\Omega_0$ requires some  
observation time $T_1$.  When $x^{(1)}$  is deeper inside $\Omega_0$
the observation time increases.  When the point  $x^{(1)}$  approaches $S_\Delta$,  i.e.  
$g_{00}(x^{(1)})\rightarrow 0$,  the needed observation time tends to infinity.  One can see this
from the fact that either forward time-like ray or backward time-like ray tends to infinity in $x_0$  when 
$x^{(1)}\rightarrow S_\Delta$  (c.f. [ER]). 

\section{Black holes and inverse problems II.}
\label{section 9}
\init

In this section we consider black or white holes inside the ergosphere.

Suppose $S_0\times\R$  is a characteristic surface,  $n\geq 2,$  and $S_0\subset S_\Delta$
where $S_\Delta$
is the ergosphere.  Suppose the condition  (\ref{eq:5.3}) on $S_0$ holds.  Consider
$v(x_0,x)$  in $\Omega_{ext}\times(0,T)$  such that 
\begin{equation}                                    \label{eq:9.1}
Lv=f,\ \ x\in\Omega_{ext}\times(0,T),
\end{equation}
\begin{equation}                                    \label{eq:9.2}
v(0,x)=\varphi_0(x),\ \ v_{x_0}(0,x)=\varphi_1(x),\ \ x\in\Omega_{ext}.
\end{equation}
We want to get an estimate of $v(x_0,x)$  in $\Omega_{ext}\times(0,T)$  in terms of 
$\varphi_0,\varphi_1,f$.  Let $\hat\varphi_0,\hat\varphi_1,\hat f$ be smooth extensions of
$\varphi_0,\varphi_1$ to $\R^n$  and of $f$  to $\R^n\times(0,T)$  such that
\begin{equation}                                     \label{eq:9.3}
\hat E(\hat\varphi_0,\hat\varphi_1)\leq 2 E(\varphi_0,\varphi_1),\ \ \ 
\|\hat f\|_{0,\R^n\times(0,T)}\leq 2\|f\|_{0,\Omega_{ext}\times(0,T)},    
\end{equation}
where   $E(\varphi_0,\varphi_1)=\int_{\Omega_{ext}}(\sum_{j=1}^n\varphi_{0x_j}^2+
\varphi_{1}^2)dx$   and $\tilde E(\tilde\varphi_0,\tilde\varphi_1)$
is a similar integral over $\R^n$.  Since $L$ is strictly hyperbolic there exists $\hat u$  in
$\R^n\times(0,T)$  such that
\begin{eqnarray}                                  \label{eq:9.4}
L\hat u=\hat f \ \ \ \mbox{in}\ \ R^n\times(0,T),
\\
\nonumber
\hat u(0,x)=\hat \varphi_0(x),\ \ \hat u_{x_0}(0,x)=\hat\varphi_1(x),\ \ x\in \R^n,
\end{eqnarray}
and
\begin{eqnarray}                                   \label{eq:9.5}
\max_{0\leq x_0\leq T}(\|\hat u(x_0,\cdot)\|_{1,\R^n}^2+\|\hat u_{x_0}(x_0,\cdot)\|_{0,\R^n}^2)
\\
\nonumber
\leq  C \hat E(\hat\varphi_0,\hat\varphi_1)
+C\int_0^T\int_{R^n}|\hat f|^2dxdx_0.
\end{eqnarray}
Replacing $v=u+\hat u$  we get that $u(x_0,x)$  satisfies
\begin{eqnarray}                                    \label{eq:9.6}
Lu=0\ \ \ \mbox{in}\ \ \Omega_{ext}\times(0,T),
\\
\nonumber
u(0,x)=u_{x_0}(0,x)=0,\ \ x\in \Omega_{ext}.
\end{eqnarray}
Therefore it remains to show that if $u(x_0,x)$  satisfies (\ref{eq:9.6})  then $u(x_0,x)=0$ in
$\Omega_{ext}\times(0,T)$.  Note that we could use the same approach in \S 8 too.

Let $T$ be small.  Denote by $\Gamma_1$  the characteristic surface 
different from $S_0\times \R$  and passing through 
$S_0\times\{x_0=T\}$.
Let $D_T$ be the domain bounded by $\Gamma_1,\ \Gamma_2=S_0\times[-\e,T]$  and
$\Gamma_3=\{x_0=-\e\}$.   For arbitrary point  $x^{(0)}\in S_0$  denote by $D_{0T}$
the intersection of $D_T$  with 
$\Sigma(x^{(0)})=\{(x_0,x): |x-x^{(0)}-(x-x^{(0)})\cdot\nu(x^{(0)})\nu(x^{(0)})|<\e\}$ , where
$\nu(x^{(0)})$  is the unit outward normal to $S_0$ at $x^{(0)}$.
Let $\alpha_j(x_0,x)\in C^\infty(D_T),\ \sum_{j=1}^N\alpha_j\equiv 1$ in $D_T,
\ \mbox{supp\ }\alpha_j\subset D_{jT}$,  
where $D_{jT}$ corresponds to $x^{(j)}\in S_0$  instead of $x^{(0)},\ \alpha_j=0$
in a neighborhood of the boundary of $\Sigma(x^{(0}))$.

Let $\alpha_0$  be any of $\alpha_j,\ 1\leq j\leq N$.
Denote $u_0=\alpha_0u$.  Then
\begin{equation}                                  \label{eq:9.7}
Lu_0=f_0,\ \ u_0(0,x)=u_{x_0}(0,x)=0,
\end{equation}
where $f_0=L'u,\ \mbox{ord\ } L'\leq 1,\ \mbox{supp\ }f_0\subset D_{0T}$.
We introduce local coordinates in a neighborhood
$B_{\e,T}=\{(x_0,x):\  x_0\in[-\e,T],|x-x^{(0)}|<2\e\}$.
Let $s=\varphi(x)$  be the solution of the eiconal equation 
\begin{equation}                                  \label{eq:9.8}
\sum_{j,k=1}^ng^{jk}(x)\varphi_{x_j}\varphi_{x_k}=0\ \ \ \mbox{in}\ \ \ B_{\e,T}.
\end{equation}
Since $S_0$  is inside the ergosphere,  $\varphi(x)$  exists when $\e$  and $T$ are
small.  We choose $s=\varphi(x)$  such that $\varphi(x)=0$  is
the equation of $S_0$  near $x^{(0)}$.

Let $\tau=\psi(x_0,x)$  be the solution of the following eiconal equation:
\begin{equation}                              \label{eq:9.9}
\sum_{j,k=0}^ng^{jk}(x)\psi_{x_j}\psi_{x_k}=0
\end{equation} 
with the initial data
\begin{equation}                           \label{eq:9.10}
\psi(x_0,x)|_{\Gamma_2}=T-x_0.
\end{equation}
Finally denote by $y_j=\varphi_j(x_0,x),\ 1\leq j\leq n-1$,
the solution of the equation
\begin{equation}                           \label{eq:9.11}
\sum_{j,k=0}^ng^{jk}(x)\psi_{x_j}\varphi_{px_k}=0\ \ \ \mbox{near}\ \ \ x^{(0)},\ \ \
1\leq p\leq n-1,
\end{equation}
with  the initial condition
\begin{equation}                           \label{eq:9.12}
\varphi_p(x_0,x)|_{\Gamma_2}=s_p(x),\ \ \ 1\leq p\leq n-1,
\end{equation}
where $s_1(x),...,s_{n-1}(x)$  are coordinates on $S_0$  near $x^{(0)},\ \frac{Dx}{D(s,s_1,...,s_{n-1})}
\neq 0 $  in $B_\e$.  Note that $\varphi_p$ does not depend on $x_0$ and $\psi(x_0,x)=T-x_0+\psi_1(x)$,
where $\psi_1(x)$  also does not depend on $x_0$.

We shall make the change of coordinates in $D_{0T}$
\begin{eqnarray}                                \label{eq:9.13}
s=\varphi(x),
\\
\nonumber
\tau=\psi(x_0,x),
\\
\nonumber
y_j=\varphi_j(x),\ \ 1\leq j\leq n-1.
\end{eqnarray}
Note that
the Jacobian $\frac{D(x_0,x)}{D(s,\tau,y')}\neq 0$  in $B_\e$  where 
$y'=(y_1,...,y_{n-1})$.  Rewrite  $Lu_0$  in $(s,\tau,y')$  coordinates (c.f. [E1], [E4]).
We get
\begin{eqnarray}                                 \label{eq:9.14}
\ \ \ \ \ \ \ \ \ \ \ \ 
\hat L\hat u_0=\frac{1}{\sqrt{|\hat g|}}\frac{\partial}{\partial s}\sqrt{|\hat g|}\hat g^{s\tau}(s,\tau,y')
\frac{\partial\hat u_0}{\partial\tau} +
\frac{1}{\sqrt{|\hat g|}}\frac{\partial}{\partial \tau}\sqrt{|\hat g|}\hat g^{s\tau}(s,\tau,y')
\frac{\partial\hat u_0}{\partial s}
\nonumber
\\
\nonumber
+\sum_{j=1}^{n-1}\frac{1}{\sqrt{|\hat g|}}\frac{\partial}{\partial s}\sqrt{|\hat g|}\hat g^{sj}(s,\tau,y')
\frac{\partial\hat u_0}{\partial y_j}+
\sum_{j=1}^n\frac{1}{\sqrt{|\hat g|}}\frac{\partial}{\partial y_j}\sqrt{|\hat g|}\hat g^{sj}(s,\tau,y')
\frac{\partial\hat u_0}{\partial s}
\\
+\sum_{j,k=1}^{n-1}\frac{1}{\sqrt{|\hat g|}}\frac{\partial}{\partial y_j}\sqrt{|\hat g|}\hat g^{jk}(s,\tau,y')
\frac{\partial\hat u_0}{\partial y_k}
\stackrel{def}{=} \hat L_1\hat u_0 +\hat L_2 \hat u_0,
\end{eqnarray}
where $L_2$  is the last sum in
(\ref{eq:9.14}) and $L_1$  are the remaining sum.  Note that $\hat g^{ss}=\hat g^{\tau\tau}=\hat g^{\tau j}=0,\ 
1\leq j\leq n-1$,  because of  (\ref{eq:9.8}), (\ref{eq:9.9}), (\ref{eq:9.11}).
In 
(\ref{eq:9.14}) $\hat u_0(s,\tau,y')=u_0(x_0,x)$  where $(x_0,x)$  and $(s,\tau,y')$  are related by
(\ref{eq:9.13}).
Since $T$ and $\e$ are small we can introduce $(s,\tau,y')$  coordinates  in $D_{0T}$.  Denote by 
$\hat D_{0T}$  the image of $D_{0T}$ in $(s,\tau,y')$  coordinates.

Let $\partial\hat \Sigma_0$  be the image of $\partial\Sigma(x^{(0)})$  in $(s,\tau,y')$  
coordinates.  Since $u_0=0$  near $\partial\Sigma(x^{(0)})$  we have that 
$\hat u_0=\hat\alpha_0\hat u=0$  near $\partial\hat \Sigma_0$.  Note also that 
$u_0=u_{0x_0}=0$  for $x_0=0$  and we extend $u_0$  by zero for $x_0<0$.
Since $\varphi_x(x^{(0)})$  is the outward normal to $\Omega_{int}$ we have
$s=\varphi(x)\geq 0$  on $\overline \Omega_{ext}$ near $S_0$.
Since 
$\tau=\psi(T,x)=0 $   on $\Gamma_2$  and $\psi_{x_0}|_{\Gamma_2}=-1$  we have that 
 $\tau\leq 0$  in 
$\hat D_{0T}$.  

Denote by $(\hat u,\hat v)$  the $L^2$  inner product in $\hat D_{0T}$.  Consider
\begin{equation}                                           \label{9.15}
(\hat L\hat u_0-\hat f_0,\hat g^{s\tau}\hat u_{0\tau}+\sum_{j=1}^{n-1}\hat g^{s0}\hat u_{0y_j}
-\hat g^{s\tau}\hat u_{0s})=0
\end{equation}
Let 
\begin{eqnarray}                                       \label{eq:9.16}
I_1=(\hat L_1\hat u_0,\hat g^{s\tau}\hat u_{0\tau}+\sum_{j=1}^{n-1}\hat g^{s0} \hat u_{0y_j})
\\
\nonumber
=\int_{\hat D_{0T}} \frac{\partial}{\partial s}(\hat g^{s\tau}\hat u_{0\tau}+
\sum_{j=1}^{n-1}\hat g^{s0}\hat u_{0y_j})^2 dsd\tau dy' +Q_1(\hat u_0,\hat u_0),
\end{eqnarray}
where  $Q_1$ has an estimate
\begin{equation}                                 \label{eq:9.17}
|Q_1(\hat u_0,\hat u_0)|\leq 
C\int_{\hat D_{0T}}(\hat u_{0s}^2+\hat u_{0\tau}^2 +\sum_{j=1}^{n-1}\hat u_{0y_j}^2)dsd\tau dy'.
\end{equation}
Therefore 
\begin{equation}                               \label{eq:9.18}
I_1=-\int_{\Gamma_2}(\hat g^{s\tau}\hat u_{0\tau}+
\sum_{j=1}^{n-1}\hat g^{s0}\hat u_{y_j})^2 d\tau dy' +Q_1(\hat u_0,\hat u_0).
\end{equation}

Denote 
\begin{equation}                                     \label{eq:9.19}
I_2=
\left
(\frac{1}{\sqrt{|\hat g|}}\frac{\partial}{\partial s}\sqrt{|\hat g|}\hat g^{s\tau}(s,\tau,y')
\frac{\partial\hat u_0}{\partial\tau} +
\frac{1}{\sqrt{|\hat g|}}\frac{\partial}{\partial \tau}\sqrt{|\hat g|}\hat g^{s\tau}(s,\tau,y')
\frac{\partial\hat u_0}{\partial s},- \hat g^{s\tau}\hat u_{0s}
\right
).
\end{equation}
We have
\begin{equation}                                       \label{eq:9.20}
I_2=-\int_{D_{0T}}\frac{\partial}{\partial\tau}\left((\hat g^{s\tau})^2\hat u_{0s}^2\right)dsd\tau dy'
+Q_2(\hat u_0,\hat u_0).
\end{equation}
Therefore 
\begin{equation}                                     \label{eq:9.21}
I_2=-\int_{\Gamma_1}(\hat g^{s\tau})^2\hat u_{0s}^2ds dy'
+Q_2(\hat u_0,\hat u_0).
\end{equation}
Denote 
\begin{equation}                                      \label{eq:9.22}
I_3=(\hat L_2\hat u_0,\hat g^{s\tau}(\hat u_{0\tau}-\hat u_{0s})).
\end{equation}
Integrating by parts in $y_j$ (note that $\hat u_0=0$  near $\partial\hat\Sigma_0$  and 
 $u_0=0$  for $x_0<0)$  we get
\begin{eqnarray}                                          \label{eq:9.23}
I_3=-\int_{\hat D_T}\sum_{j,k=1}^{n-1}\hat g^{s\tau}\hat g^{jk}(\hat u_{0sy_j}
-\hat u_{0\tau y_j})\hat u_{0y_k}dsd\tau dy'
+Q_3
\\
\nonumber
=\frac{1}{2}\int_{\Gamma_1}(\sum_{j,k=1}^{n-1}g^{jk}\hat u_{0y_j}\hat u_{0y_k})
\hat g^{s\tau}dsdy' 
+\frac{1}{2}\int_{\Gamma_2}(\sum\hat g^{jk}\hat u_{0y_j}\hat u_{0y_k})
\hat g^{s\tau}d\tau dy'
+Q_4.
\end{eqnarray}
Note that 
\begin{equation}                                \label{eq:9.24}
I_4=(\hat L_2 \hat u_0,\sum_{j=1}^{n-1}\hat g^{sj}\hat u_{0y_j})=Q_5(\hat u,\hat u)
\end{equation}
since it can be represented as a divergence 
of quadratic form in $\hat u_{0y_j}$ (c.f. (\ref{eq:8.15}).  Also 
$$
|(\hat f_0,\hat g^{s\tau}(\hat u_{0\tau}-\hat u_{0s}) +\sum_{j=1}^{n-1}\hat g^{js}\hat u_{0y_j})|
\leq C\int_{\hat D_{0T}}|\hat f_0|^2dsd\tau dy' +Q_6(\hat u_0,\hat u_0).
$$
Therefore we have
\begin{eqnarray}                                         \label{eq:9.25}
\frac{1}{2}\int_{\Gamma_2}
(-\sum_{j,k=1}^{n-1}\hat g^{jk}\hat u_{0y_j}\hat u_{0y_k}\hat g^{s\tau}
)d\tau dy'
+
\frac{1}{2}\int_{\Gamma_1}
(-\sum_{j,k=1}^{n-1}\hat g^{jk}\hat u_{0y_j}\hat u_{0y_k}\hat g^{s\tau}
)ds dy'
\nonumber
\\
\nonumber
+\int_{\Gamma_2}
( \hat g^{s\tau}\hat u_{0\tau}+\sum_{j=1}^{n-1}\hat g^{sj}\hat u_{0y_j} 
)^2d\tau dy'
+\int_{\Gamma_1}(\hat g^{s\tau})^2\hat u_{0s}^2dsdy'
\\
\leq Q_7(\hat u_0,\hat u_0)+C\int_{\hat D_{0T}}|\hat f_0|^2dsd\tau dy'
\end{eqnarray}
Note that $\hat g^{s\tau}>0$ and
\begin{equation}                                   \label{eq:9.26}                 
-\sum_{j,k=1}^{n-1}\hat g^{jk}\hat u_{0y_j}\hat u_{0y_k}
\geq 
C\sum_{j=1}^{n-1}\hat u_{0y_j}^2 \ \ \mbox{in}\ \ \hat{D}_{0T}.
\end{equation}
Therefore
\begin{equation}                                  \label{eq:9.27}
\int_{\Gamma_1}
[(\hat g^{s\tau})^2\hat u_{0s}^2 
-\frac{1}{2}\sum_{j,k=1}^{n-1}\hat g^{jk}\hat u_{0y_j}\hat u_{0y_k} \hat g^{s\tau}
]dsdy'
\end{equation}
is equivalent to 
$\int_{\Gamma_1}(\hat u_{os}^2+\sum_{j=1}^{n-1}\hat u_{0y'}^2)dsdy'$,  i.e. it is 
equivalent to the norm $\|\hat u_0\|_{1,\Gamma_1}^2$  in $H^1(\Gamma_1)$.   Analogously
\begin{equation}                                             \label{eq:9.28}
\int_{\Gamma_2}
[
( \hat g^{s\tau}\hat u_{0\tau}+
\sum_{j=1}^{n-1}\hat g^{sj}\hat u_{0y_j}
)^2
-\frac{1}{2}\sum_{j,k=1}^{n-1}\hat g^{jk}\hat u_{0y_j}\hat u_{0y_k}\hat g^{s\tau}
]d\tau dy'
\end{equation}
is equivalent to the norm
$$
\|u_0\|_{1,\Gamma_2}^2
=
\int_{\Gamma_2}
( \hat u_{0\tau}^2+
\sum_{j=1}^{n-1}\hat u_{0y_j}^2)d\tau dy'
$$
  in $H^1(\Gamma_2)$.
Therefore (\ref{eq:9.25})  is equivalent to
\begin{equation}                           \label{eq:9.29}
\|\hat u_0\|_{1,\Gamma_1}^2+\|\hat u_0\|_{1,\Gamma_2}^2\leq Q_8(\hat u_0,\hat u_0)+
C\int_{\hat D_{0T}}|\hat f_0|^2dsd\tau dy'.
\end{equation}
Denote by $D_{0T,t}$  the intersection of $D_{0T}$  with the half-space $x_0\geq t$.
Integrating by parts in the integral
\begin{equation}                               \label{eq:9.30}
0=\int_{D_{0T,t}}(Lu_0-f_0)(g^{00}u_{0x_0}+\sum_{j=1}^ng^{0j}(x)u_{0x_j})dxdx_0
\end{equation}
we obtain (c.f. [E1], [E3])
\begin{eqnarray}                              \label{eq:9.31}
\int_{D_{0T}\cap\{x_0=t\}}
[(\sum_{j=0}^ng^{j0}u_{0x_j}(t,x))^2\ \ \ \ \ \ \ \ \ \ \ \ \ \ \ \ \ \ \ \
\ \ \ \ \ \ \ \ \ \ 
\\
+ \ \ (\sum_{j=1}^ng^{j0}u_{0x_j}(t,x))^2
-
\sum_{j,k=1}^ng^{jk}u_{0x_j}(t,x)u_{0x_k}(t,x)]dx 
\nonumber
\\
\leq C(\|u_0\|_{1,\Gamma_{1t}}^2+\|u_0\|_{1,\Gamma_{2t}}^2)
+C\int_{D_{0T,t}}(\sum_{j=0}^nu_{0x_j}^2)dxdx_0 
+C\int_{D_{0T,t}}|f_0|^2dxdx_0,
\nonumber
\end{eqnarray}
where
$\Gamma_{jt}$  is the intersection of $\Gamma_j$  with $x_0\geq t,\ j=1,2.$
Note that the integral in the left hand side of (\ref{eq:9.31}) is equivalent to
\begin{equation}                                   \label{eq:9.32}
\int_{D_{0T}\cap\{x_0=t\}}
(\sum_{j=0}^nu_{0x_j}^2(t,x))dx.
\end{equation}
Rewriting (\ref{eq:9.29})  in $(x_0,x)$  coordinates and combining with (\ref{eq:9.31})
we get
\begin{equation}                                 \label{eq:9.33}
\max_{0\leq t\leq T}\int_{D_{0T}\cap\{x_0=t\}}
(\sum_{j=0}^nu_{0x_j}^2)dx
\leq
C\int_{D_{0T}}
(\sum_{j=0}^nu_{0x_j}^2)dx_0dx+C\int_{D_{0T}}|f_0|^2dx_0dx.
\end{equation}
Let $\{\alpha_j(x)\}_{j=1,...,N}$  be as above.  Denote $u_j=\alpha_j u$.
Applying (\ref{eq:9.33}) with $u_j=\alpha_ju$ instead of $u_0=\alpha_0 u$  
and using that $\sum_{j=1}^N\alpha_j=1$  in $D_T$ we get
\begin{eqnarray}                                 \label{eq:9.34}
\max_{0\leq t\leq T}\int_{D_{T}\cap\{x_0=t\}}
(\sum_{j=0}^nu_{x_j}^2(t,x))dx
\ \ \ \ \ \ \ \ \ \ \ \ \ \ \ \ \ \ \ 
\\
\leq
C\int_{D_{T}}
(\sum_{j=0}^nu_{x_j}^2(x_0,x))dx_0dx+C\sum_{j=1}^N\int_{D_{jT}}|f_j|^2dx_0dx.
\nonumber
\end{eqnarray}
where
$f_j=(L\alpha_j-\alpha_j L)u$,
\begin{equation}                                    \label{eq:9.35}
|f_j|\leq C\sum_{j=0}^n|u_{x_j}|.
\end{equation}
Therefore 
\begin{eqnarray}                                 \label{eq:9.36}
\max_{0\leq t\leq T}\int_{D_{T}\cap\{x_0=t\}}
(\sum_{j=0}^nu_{x_j}^2(t,x))dx
\ \ \ \ \ \ \ \ \ \ \ \ \ \ \ \ \ \ \ 
\\
\leq
C\int_{D_{T}}
(\sum_{j=0}^nu_{x_j}^2(x_0,x))dx_0dx \leq C T
\max_{0\leq t\leq T}\int_{D_{T}\cap\{x_0=t\}}
(\sum_{j=0}^nu_{x_j}^2(t,x))dx
\nonumber
\end{eqnarray}
Since $T$  is small we conclude that $u=0$  in $D_T$.  Take any $T_1<T$.  Then there exists 
$\delta_1>0$  such that $S_{0\delta_1}\times[0,T_1]\subset D_T$  where $S_{0\delta_1}$
is a $\delta_1$-neighborhood of $S_0$.  Therefore $u(x_0,x)=0$   in $S_{0\delta_1}\times[0,T_1]$  
and we can  extend $u(x_0,x)$  by zero in $\Omega_{int}\times[0,T_1]$.  Then $Lu=0$
in $\R^n\times(0,T_1)$  and $u(0,x)=u_{x_0}(0,x)=0$  in $\R^n$.  By the uniqueness of the hyperbolic 
Cauchy problem (c.f. (\ref{eq:9.5})) we have $u=0$  in $\R^n\times (0,T_1)$.  Repeating the same 
arguments on
$(T_1,2T_1)$, etc.,  we get that $u=0$  in $\Omega_{ext}\times(0,T)$  for any $T>0$.
Therefore $v=\hat u$  in $\Omega_{ext}\times (0,T)$  where $v(x_0,x)$  satisfies 
(\ref{eq:9.1}),  (\ref{eq:9.2}) and
$\hat u$  satisfies (\ref{eq:9.4})  in $\R^n\times (0,T)$.  Then  (\ref{eq:9.5}) implies that
\begin{eqnarray}                                \label{eq:9.37}
\max_{0\leq x_0\leq T}(\|v(x_0,\cdot)\|_{1,\Omega_{ext}}^2  + \|v_{x_0}(x_0,\cdot)\|_{0,\Omega_{ext}}^2)
\\
\nonumber
\leq C(\|\varphi_0\|_{1,\Omega_{ext}}^2+\|\varphi_1\|_{0,\Omega_{ext}}^2 
+\int_0^T\|f_(x_0,\cdot)\|_{0,\Omega_{ext}}^2dx_0).
\end{eqnarray}
Therefore we proved an analogue of Theorem \ref{theo:8.3}:
\begin{theorem}                                 \label{theo:9.1}
Let $S_0$  be a characteristic surface inside the ergosphere $S_\Delta$   and let $v(x_0,x)$
satisfies (\ref{eq:9.1}), (\ref{eq:9.2}).  Suppose  (\ref{eq:5.3})
holds.  Then $v(x_0,x)$ satisfies (\ref{eq:9.37}).
\end{theorem}
Note that (\ref{eq:9.37})
implies that $D_+(\Omega_{int}\times \R)\subset\overline{\Omega}_{int}\times \R$,  i.e.  that 
$\Omega_{int}\times \R$  is a black hole.

As in the case of Theorem \ref{theo:8.1}  we can take into account boundary condition on 
$\partial\Omega_0$ and prove that estimate of the form 

\begin{eqnarray}                                    \label{eq:9.38}
\max_{0\leq x_0\leq T}(\|v(x_0,\cdot)\|_{1,\Omega_{ext}\cap\Omega}^2+
\|v_{x_0}(x_0,\cdot)\|_{0,\Omega_{ext}\cap\Omega}^2)
\\
\nonumber
\leq C\int_0^T\|f(x_0,\cdot)\|_{0,\Omega_{ext}\cap\Omega}^2dx_0
+[g]_{1,\partial\Omega_0\times(0,T)}^2.
\end{eqnarray}
When (\ref{eq:5.2}) holds
we get that the domain $D_T$  will be contained in $\Omega_{int}\times(0,T)$.
In this case the proof similar to the proof of Theorem \ref{theo:9.1}  gives an estimate of the form 
(\ref{eq:9.38})
in $(\Omega_{int}\cap\Omega)\times(0,T)$,
  i.e.  in this case $\Omega_{int}\times\R$ is a white hole.


\begin{thebibliography}{9999}
\bibitem[B]{} Belishev, M., 2007, Recent progress in the boundary control method,
Inverse Problems 23, No 5, R1-67
\bibitem[B1]{} Belishev, M., 1997, Boundary control in reconstruction
of manifolds and metrics (the BC method),
Inverse Problems 13, R1-R45
\bibitem[CH]{} Courant, R., Hilbert, D., Methods of Mathematical Physics,
vol. II (1962),  New York, London  
\bibitem[E1]{} Eskin, G., 2006, A new approach to the hyperbolic
inverse problems,  Inverse problems, vol. 22, No. 3, 815-831
\bibitem[E2]{} Eskin, G., 2007, A new approach to the hyperbolic
inverse problems II: global step,  
 Inverse Problems 23, 2343-2356
\bibitem[E3]{} Eskin, G., 2007, Inverse  hyperbolic problems 
with time-dependent coefficients,   Comm. in PDE 32, 1737-1758
\bibitem[E4]{} Eskin, G., 2008, 
Optical Aharonov-Bohm effect:  inverse hyperbolic problem approach,
Comm. Math. Phys., 284,  317-343 
\bibitem[E5]{} Eskin, G., 2008,Iverse hyperbolic problems and optical black holes,
ArXiv:0809.3987
\bibitem[E6]{} Eskin, G., 1987, Mixed initial-boundary value prolems for second
order hyperbolic equations, Comm. in PDE, 12:503-587
\bibitem[ER]{} Eskin, G., Ralston, J., 2009, On the determination
boundaries for hyperbolic equations,  ArXiv:0902.4497  
\bibitem[G]{} Gordon, W., 1923, Ann. Phys.  (Leipzig) 72, 421
\bibitem[H]{} Hartman, F., 1964, Ordinary differential equations (New York, J.Wiley \& son)
\bibitem[Ho]{} Hormander, L., 1985, The Analysis of Linear Partial 
Differential Operators III (Berlin: Springer)
\bibitem[KKL]{} Katchalov, A., Kurylev, Y., Lassas, M., 2001,
Inverse boundary spectral problems (Boca Baton : Chapman\&Hall)
 \bibitem[LP]{} Leonhardt, V., Piwnicki, P., 1999,  Phys.  Rev. A60,  4301
\bibitem[NVV]{} Novello, M., Visser, M.,  Volovik, G. (editors),
Artificial black holes, 2002,  World Scientific,  Singapore.
\bibitem[T]{} Tataru, D., 1995, Unique continuation for solutions
to PDE, Comm. in PDE 20, 855-884
\bibitem[V]{} Visser, M., 1998,  Acoustic black holes, horizons,  ergospheres and Hawking
radiation,  Classical quantum gravity 15, No. 6, 1767-1791. 














\end{thebibliography}
\end{document}